\documentstyle[11pt]{article}
\input amssym.def
\input amssym.tex
\input psfig
\makeatletter
\def\@sect#1#2#3#4#5#6[#7]#8{\ifnum #2>\c@secnumdepth
     \def\@svsec{}\else 
     \refstepcounter{#1}\edef\@svsec{\csname the#1\endcsname.\hskip .75em }\fi
     \@tempskipa #5\relax
      \ifdim \@tempskipa>\z@ 
        \begingroup #6\relax
          \@hangfrom{\hskip #3\relax\@svsec}{\interlinepenalty \@M #8\par}%
        \endgroup
       \csname #1mark\endcsname{#7}\addcontentsline
         {toc}{#1}{\ifnum #2>\c@secnumdepth \else
                      \protect\numberline{\csname the#1\endcsname}\fi
                    #7}\else
        \def\@svsechd{#6\hskip #3\@svsec #8\csname #1mark\endcsname
                      {#7}\addcontentsline
                           {toc}{#1}{\ifnum #2>\c@secnumdepth \else
                             \protect\numberline{\csname the#1\endcsname}\fi
                       #7}}\fi
     \@xsect{#5}}

\def\section{\@startsection {section}{1}{\z@}{-3.5ex plus -1ex minus 
 -.2ex}{2.3ex plus .2ex}{\normalsize\bf}}
\makeatother

\def\binom#1#2{{#1}\choose{#2}}
\newcommand{\sA}{{\cal A}}
\newcommand{\sB}{{\cal B}}
\newcommand{\sD}{{\cal D}}
\newcommand{\sL}{{\cal L}}
\newcommand{\af}{\alpha}
\newcommand{\dd}{\ldots}
\newcommand{\eeq}{\end{equation}}
\newcommand{\beql}[1]{\begin{equation}\label{#1}}
\newcommand{\eqn}[1]{(\ref{#1})}

\newcommand{\hsp}{\hspace*{\parindent}}

\begin{document}
\begin{center}
{\Large {\bf Some Canonical Sequences of Integers}} \\
\vspace{1\baselineskip}
{\em M. Bernstein} and {\em N. J. A. Sloane} \\
\vspace{.25\baselineskip}
Mathematical Sciences Research Center \\
AT\&T Bell Laboratories \\
Murray Hill, NJ 07974 USA \\
\vspace{2\baselineskip}
{\sc Dedicated to Professor J.~J. Seidel} \\
\vspace{2\baselineskip}
{\bf Abstract} \\
\vspace{.5\baselineskip}
\end{center}
\setlength{\baselineskip}{1.5\baselineskip}

Extending earlier work of R. Donaghey and P. J. Cameron, we investigate
some canonical ``eigen-sequences'' associated with transformations of integer sequences.
Several known sequences appear in a new setting:
for instance the sequences (such as 1, 3, 11, 49, 257, 1531, $\dd$) studied by
T. Tsuzuku, H.~O. Foulkes and A. Kerber in connection with
multiply transitive groups are eigen-sequences for the binomial transform.
Many interesting new sequences also arise, such as 1, 1, 2, 26, 152, 1144, $\dd$,
which shifts one place left when transformed by the Stirling numbers of the second kind,
and whose exponential generating function satisfies $\sA ' (x) = \sA (e^x -1) +1$.
\clearpage
\large\normalsize
\renewcommand{\baselinestretch}{1}
\thispagestyle{empty}
\setcounter{page}{1}
\begin{center}
{\Large {\bf Some Canonical Sequences of Integers}} \\
\vspace{1\baselineskip}
{\em M. Bernstein} and {\em N. J. A. Sloane} \\
\vspace{.25\baselineskip}
Mathematical Sciences Research Center \\
AT\&T Bell Laboratories \\
Murray Hill, NJ 07974 USA \\
\vspace{1.5\baselineskip}
\end{center}
\setlength{\baselineskip}{1.5\baselineskip}

\section{Transformations of sequences}
\hsp
The purpose of this note is to present certain integer sequences that are associated
in a canonical way with various sequence transformations.
This is an extension of earlier work of R. Donaghey \cite{Dona76}, who was concerned with the transformation denoted by BINOMIAL in the present paper, and of P.~J. Cameron \cite{Came89},
who mainly considered the transforms that we call EULER, INVERT and WEIGH.
Cameron says that the last of these transformations
was suggested to him by Professor Seidel.
We hope that Professor Seidel will enjoy the sequences presented here.

A number of papers have appeared in recent years in the combinatorial literature dealing with generalizations of the Catalan and Bell numbers (e.g. \cite{Dona77}, \cite{DoSh77}, \cite{Roge77}, \cite{StWa79}).
The present paper provides another way to view certain of these sequences, as well as suggesting still further generalizations.

We shall not give any proofs of the identification of certain of our sequences with those in the literature --- in most cases these are easily obtained from known properties of the sequences.
The majority of the sequences discussed here are not in \cite{HIS}, although they will appear in \cite{SlPl94}.

A typical sequence will be written $a= [a_0, a_1 , \dd ]$, with ordinary
generating function (o.g.f.) $A (x) = \sum\limits_{n=0}^\infty a_n x^n$ and exponential generating function (e.g.f.) $\sA (x) = \sum\limits_{n=0}^\infty a_n x^n / n!$.
The logarithmic derivative $\sA ' (x) / \sA (x)$ is denoted by $\sL (x)$.
Unless indicated otherwise, ``sequence'' means ``sequence of real numbers''.
The following list contains the transformations that we shall consider.
Here $a$ is the original sequence and $b$ the transformed sequence.
In each case $b$ is integral if $a$ is.

The initial transformations map a sequence $a = [a_0 , a_1 , \dd ]$ to $b= [b_0 , b_1 , \dd ]$.
\paragraph{BINOMIAL}
\beql{eq1}
b_n = \sum_{k=0}^n {\binom{n}{k}} a_k , ~~a_n = \sum_{k=0}^n
(-1)^{n-k} {\binom{n}{k}} b_k ~,
\eeq
with two equivalent forms
\beql{eq2}
B(x) = \frac{1}{1-x} A \left( \frac{x}{1-x} \right) ~,~~~
\sB (x) = e^x \sA (x) ~.
\eeq
References: \cite[p.~4]{RCI}, \cite[p~192]{GKP},
\cite[p.~6]{GK90}, \cite{Dona76}, \cite{DoSh77} (in the latter reference the inverse of this transformation appears under the name ``Euler transformation'').
If $a_n$ enumerates a class of structures on $n$ points in which every element is involved,
$b_n$ enumerates the same type of structure when only a subset
of the points needs to be involved.

This transformation also arises when analyzing sequences by means of difference
tables.
The familiar difference table associated with a sequence
$[c_0 , c_1 , \dd ]$ is its difference table of depth 1.
By using the leading diagonal of this table as the top row of a new table,
we obtain the difference table of depth 2, and so on.

For example the difference tables of depths 1 and 2 for $[1, 3, 9, 27, \dd ]$ are
$$
\begin{array}{ccccccc}
1 & ~ & 3 & ~ & 9 & ~ & 27 \\
~ & 2 & ~ & 6 & ~ & 18 & ~ \\
~ & ~ & 4 & ~ & 12 & ~ \\
~ & ~ & ~ & 8 
\end{array}
~~~~~
~~
\begin{array}{ccccccc}
1 & ~ & 2 & ~ & 4 & ~ & 8 \\
~ & 1 & ~ & 2 & ~ & 4 \\
~ & ~ & 1 & ~ & 2 \\
~ & ~ & ~ & 1
\end{array}
$$
To recover the sequence from the leading diagonal of the difference table of depth $r$, we apply the binomial transform $r$ times.
In the example we go from $[1,1,1,1, \dd ]$ to $[1,2,4,8, \dd]$ to
$[1,3,9,27, \dd ]$.
\paragraph{STIRLING}
$$b_n = \sum_{k=0}^n S(n,k) a_k ,~~
a_n = \sum_{k=0}^n s(n,k) b_k ~,
$$
or equivalently
\beql{eq3}
\sB (x) = \sA (e^x -1) ~,
\eeq
where $S(n,k)$ and $s(n,k)$ are respectively Stirling numbers of the second and first kinds (the numbers of the second kind being the more fundamental)
\cite[pp.~90,202]{RCI}
\cite[p.~296]{GKP},
\cite[p.~7]{GK90}.
If $a_n$ is the number of objects in some class with points labeled
$1,2, \dd, n$ (with all labels distinct, i.e. ordinary
labeled structures),
then $b_n$ is the number of objects with points labeled $1,2, \dd, n$ (with repetitions allowed).

The following transformations are related to the operation of convolution.
They do not in general have inverses.
\paragraph{CONV}
$$
b_n =
\sum\limits_{k=0}^n a_k a_{n-k} ,~~
B(x) = A(x)^2 ~.
$$
\paragraph{EXP-CONV}
$$b_n = \sum\limits_{k=0}^n {\binom{n}{k}} a_k a_{n-k} ,~~
\sB (x) = \sA (x)^2 ~.
$$
\paragraph{F-CONV}
$$b_n = \sum\limits_{k=0}^n F(a_k, a_{n-k})~,$$
where $F(r,s)$ is
one of LCM, GCD, AND, OR, XOR (in the last three $F$ is applied to the binary representation of its arguments).
The F-CONV transforms were suggested to us by M. Le~Brun \cite{lebr},
and will be applied only to sequences of nonnegative integers.

The remaining transformations map
$a= [a_1 , a_2, \dd ]$ to $b= [b_1 , b_2 , \dd]$.
\paragraph{M\"{O}BIUS}
$$b_n = \sum_{d|n} \mu \left( \frac{n}{d} \right) a_d ,~~
a_n = \sum_{d|n} b_d ~,$$
where $\mu$ is the M\"{o}bius function.
Two equivalent forms are
\begin{eqnarray}
\label{eq50}
\sum_{n=1}^\infty \frac{a_n}{n^s} & = & \zeta (s)
\sum_{n=1}^\infty \frac{b_n}{n^s} ~, \\
\sum_{n=1}^\infty a_n x^n & = & \sum_{n=1}^\infty
\frac{b_n x^n}{1-x^n} ~, \nonumber
\end{eqnarray}
the latter being a Lambert series \cite[p.~257]{HW1}.
\paragraph{WEIGH}
$$1+ \sum_{n=1}^\infty b_n x^n =
\prod_{n=1}^\infty (1+x^n)^{a_n} ~.$$
$b_n$ is the number of ways of getting a weight of $n$ grams
if we have $a_i$ weights of $i$ grams, $i \ge 1$
\cite[p.~174]{PoSz72}.
Cameron \cite{Came89} gives other interpretations.
\paragraph{EULER}
$$1+ \sum_{n=1}^\infty b_n x^n = \prod_{n=1}^\infty
\frac{1}{(1-x^n)^{a_n}}~.$$
If $a_n$ enumerates a class of connected structures on $n$ unlabeled nodes,
$b_n$ enumerates the same structures when connectedness is not required.
$b_n$ is also the number of ways of partitioning the integer $n$,
given that there are $a_k$ possible types of parts of size $k$,
for $k=1,2, \dd$.
Alternatively, if $a_n$ is the number of generators of degree $n$ of a graded polynomial algebra, then $b_n$ is the dimension of the $n$-th homogeneous
component, i.e. the number of linearly independent monomials of degree $n$.
WEIGH has a similar interpretation in terms of graded exterior algebras.
References \cite{BeGo71}, \cite{Rota75} and \cite{Came89} give other interpretations.
Calculations are facilitated by setting the left side equal to
$\exp \sum\limits_{n=1}^\infty c_n x^n /n$, so that
$c_n = nb_n - \sum\limits_{k=1}^{n-1} c_k b_{n-k}$, $a_n = \frac{1}{n} \sum\limits_{d|n}
\mu \left( \frac{n}{d} \right) c_d$.
\paragraph{PARTITION.}
For $1 \le a_1 \le a_2 \le \dd$,
form $c= [c_1 , c_2 , \dd]$ by discarding all duplicate entries from $a$,
so that $c$ is strictly increasing;
then
$$1+ \sum_{n=1}^\infty b_n x^n = \prod_{n=1}^\infty \frac{1}{1-x^{c_n}} ~.
$$
$b_n$ is the number of ways of partitioning $n$ into parts that can have sizes
$a_1, a_2 , \dd$~.
\paragraph{INVERT}
$$1 + \sum_{n=1}^\infty b_n x^n =
\frac{1}{1- \sum\limits_{n=1}^\infty a_n x^n} ~.$$
$b_n$ is the number of ordered arrangements of postage stamps of total value $n$
that can be formed if we have $a_i$ types of stamps of value $i$,
$i \ge 1$ \cite[p.~174]{PoSz72}.
\cite{Came89} gives several other interpretations.
\paragraph{REVERT.}
If $y= \sum\limits_{n=1}^\infty a_n x^n$ with $a_1 =1$, define $b$ by
$x = \sum\limits_{n=1}^\infty (-1)^{n+1} b_n y^n$ (so $b_1 =1$)
\cite[p.~149]{RCI}, \cite[p.~16]{AS1}, \cite{Knut92}.
\paragraph{EXP}
$$
1+ \sB (x) = \exp \sA (x) ,~~
\sA (x) = \log (1+ \sB (x)) ~.$$
This transformation relates the number connected structures on $n$ labeled
nodes to the total number of structures
\cite{BeGo71}, \cite{Rota75}.
\section{Eigen-sequences}
\hsp
Let $T$ denote any of the above transformations.
We are interested in {\em eigen-sequences} $a$ such that $T \circ a =a$, or (since in
most cases that equation has only trivial solutions)
$$
\begin{array}{rllrll}
T \circ a & = & L \circ a ,~ &
T^2 \circ a & = & L \circ a ,~~~T \circ a ~=~ L^2 \circ a ~, \\ [+.15in]
T \circ a & = & N \circ a , ~ &
T \circ a & = & M \circ a ~,
\end{array}
$$
etc., where the auxiliary operators $R$, $L$, $M$, $N$ are defined by
\begin{eqnarray*}
R \circ [a_0 , a_1 , a_2 , \dd ] & = & [1, a_0 , a_1 , \dd ] ~, \\
L \circ [a_0 , a_1 , a_2 , \dd ] & = & [a_1 , a_2 , a_3 , \dd ] ~, \\
N \circ [a_0 , a_1 , a_2 , \dd ] & = & [a_0, -a_1, -a_2, \dd ] ~, \\
M \circ [a_0 , a_1 , a_2 , \dd ] & = & [a_1, 2a_1, 2a_2 , \dd ] ~.
\end{eqnarray*}
Note that the e.g.f.'s for $R \circ a$ and $L \circ a$ are respectively
$1 + \int_0^x \sA (x) dx$ and $\sA ' (x)$.
The equations $T \circ a = L \circ a$ and $T \circ a = M \circ a$ were studied by
Cameron \cite{Came89} for the transformations EULER, INVERT and WEIGH.

Suppose $T$ is linear, in the sense that it maps $a$ to $b$ with
$$b_n = \sum_{k=0}^n D_{n,k} a_k ~.$$
Let $D$ denote the array $[D_{n,k} : n=0, \dd ; k=0, \dd ]$,
and set $D^r = [D_{n,k}^{(r)} ]$.
Then a sequence satisfying
$T^r \circ a = L^s \circ a$ has the recurrence
\beql{eq21}
a_{n+s} = \sum_{k=0}^n D_{n,k}^{(r)} a_k ~,
\eeq
and, provided $D_{n,n} =1$ for all $n$,
sequences satisfying $T \circ a =N \circ a$ or $T \circ a = M \circ a$ have recurrences
\begin{eqnarray}
\label{eq22}
a_n & = & - \frac{1}{2} \sum_{k=0}^{n-1} D_{n,k} a_k ~, \\
\noalign{or}
\label{eq23}
a_n & = & \sum_{k=0}^{n-1} D_{n,k} a_k ~,
\end{eqnarray}
respectively.

Tables~I and II show the most interesting of these eigen-sequences.
The symbol in the last column of Table~I indicates the sense in which the sequence is unique:

\noindent
$\af$:
unique sequence fixed by the operator;
furthermore any sequence converges to this one under repeated application of the operator

\noindent
$\beta$:
unique sequence beginning $1, \dd$ that is fixed by the operator

\noindent
$\gamma$:
unique sequence $[0,1, a_2 , a_3 , \dd ]$ of nonnegative integers that is mapped
to $[0, a_2 , a_3 , \dd ]$ under XOR-CONV

\noindent
$\delta$:
every sequence converges to one of these five (see below for further details)

\noindent
$\epsilon$:
lexicographically earliest integer sequence $[1= a_1 < a_2 < \cdots ]$ fixed by REVERT.
\begin{table}[htb]
\begin{center}
{\bf Table~I(a)~~~Integer eigensequences} \\
\end{center}

$$
\begin{array}{llll}
\multicolumn{1}{c}{\#} & \multicolumn{1}{c}{\mbox{Sequence}} & \multicolumn{1}{c}{\mbox{Operator}} & \multicolumn{1}{c}{\mbox{Property}} \\ [+.2in] \\
S1 & [1,1,2,5,15,52,203,877,4140,...] & R \circ \mbox{BINOMIAL} & \alpha \\
S2 & [1,1,3,11,49,257,1539,10299,...] & R \circ \mbox{BINOMIAL}^2 & \alpha \\
S3 & [1,1,4,19,109,742,5815,51193,...] & R \circ \mbox{BINOMIAL}^3 & \alpha \\
S4 & [1,1,5,29,201,1657,15821,...] & R \circ \mbox{BINOMIAL}^4 & \alpha \\
S5 & [1,1,1,2,4,9,23,65,199,654,...] & R^2 \circ \mbox{BINOMIAL} & \alpha \\
S6 & [1,1,3,13,75,541,4683,47293,...] & M^{-1} \circ \mbox{BINOMIAL} & \beta \\
\\
S7 & [1,1,2,6,26,152,1144,10742,...] & R \circ \mbox{STIRLING} & \alpha \\
S8 & [1,1,3,14,97,934,11814,188650,...] & R \circ \mbox{STIRLING}^2 & \alpha \\
S9 & [1,1,1,2,5,16,66,343,2167,...] & R^2 \circ \mbox{STIRLING} & \alpha \\
S10 & [1,1,4,32,436,9012,262760,...] & M^{-1} \circ \mbox{STIRLING} & \beta \\
\\
S11 & [1,1,2,5,14,42,132,429,1430,...] & R \circ \mbox{CONV} & \alpha \\
S12 & [1,1,4,22,140,969,7084,53820,...] & R \circ \mbox{CONV}^2 & \alpha \\
S13 & [1,1,8,92,1240,18278,285384,...] & R \circ \mbox{CONV}^3 & \alpha \\
S14 & [1,1,1,2,3,6,11,22,44,90,187,...] & R^2 \circ \mbox{CONV} & \alpha \\
\\
S15 & [1,1,2,6,24,120,720,5040,...] & R \circ \mbox{EXP-CONV} & \alpha \\
S16 & [1,1,4,28,280,3640,58240,...] & R \circ \mbox{EXP-CONV}^2 & \alpha \\
S17 & [1,1,1,2,4,10,30,100,380,1600,.. & R^2 \circ \mbox{EXP-CONV} & \alpha \\
\\
S18 & [1,1,2,5,14,40,128,369,1214,...] & R \circ \mbox{LCM-CONV} & \alpha \\
S19 & [1,1,2,3,4,6,6,11,10,18,16,...] & R \circ \mbox{GCD-CONV} & \alpha \\
S20 & [1,1,2,1,2,4,0,5,2,4,0,10,0,...] & R \circ \mbox{AND-CONV} & \alpha \\
S21 & [1,1,2,7,20,58,174,519,1550,...] & R \circ \mbox{OR-CONV} & \alpha \\
S22 & [0,1,2,4,14,38,118,338,1006,...] & R \circ \mbox{XOR-CONV} & \gamma \\
\\
S23 & [1,1,2,3,5,6,10,11,16,19,26,...] & R \circ \mbox{M\"{O}BIUS}^{-1} & \alpha \\
S24 & [1,1,3,5,10,12,24,26,43,52,78,...] & R \circ \mbox{M\"{O}BIUS}^{-2} & \alpha \\
S25 & [1,1,1,2,2,4,3,7,4,11,6,15,7,...] & R^2 \circ \mbox{M\"{O}BIUS}^{-1} & \alpha \\
S26 & [1,1,1,2,1,3,1,4,2,3,1,8,1,3,3,...] & M^{-1} \circ \mbox{M\"{O}BIUS}^{-1} & \beta
\end{array}
$$
\end{table}

\begin{table}[htb]
\begin{center}
{\bf Table~I(b)~~~Integer eigensequences (cont.)} \\
\end{center}

$$
\begin{array}{llll}
\multicolumn{1}{c}{\#} & \multicolumn{1}{c}{\mbox{Sequence}} & \multicolumn{1}{c}{\mbox{Operator}} & \multicolumn{1}{c}{\mbox{Property}} \\ [+.2in] \\
S27 & [1,1,1,2,3,6,12,25,52,113,247,...] & R \circ \mbox{WEIGH} & \alpha \\
S28 & [1,1,1,3,6,16,43,120,339,985,...] & R \circ \mbox{WEIGH}^2 & \alpha \\
S29 & [1,1,1,1,2,2,4,6,10,17,29,51,...] & R^2 \circ \mbox{WEIGH} & \alpha \\
\\
S30 & [1,1,2,4,9,20,48,115,286,719,...] & R \circ \mbox{EULER} & \alpha \\
S31 & [1,1,3,8,25,77,258,871,3049,...] & R \circ \mbox{EULER}^2 & \alpha \\
S32 & [1,1,1,2,3,6,10,20,36,72,137,...] & R^2 \circ \mbox{EULER} & \alpha \\
S33 & [1,1,2,5,12,33,90,261,766,2312,...] & M^{-1} \circ \mbox{EULER} & \beta \\
\\
S34 & [1,2,2,4,5,7,9,12,16,20,25,32,...] & \mbox{PARTITION} & \delta \\
S35 & [1,2,2,4,4,6,7,11,12,16,18,25,...] & \mbox{PARTITION}^2 & \delta \\
S36 & [1,2,2,4,4,7,8,12,13,18,21,29,...] & \mbox{PARTITION}^2 & \delta \\
S37 & [1,2,3,4,6,9,11,15,19,25,31,41,...] & \mbox{PARTITION}^2 & \delta \\
S38 & [1,2,3,5,6,10,12,17,22,29,36,48,...] & \mbox{PARTITION}^2 & \delta \\
\\
S11 & [1,1,2,5,14,42,132,429,1430,...] & R \circ \mbox{INVERT} & \alpha \\
S39 & [1,1,3,11,45,197,903,4279,20793,...] & R \circ \mbox{INVERT}^2 & \alpha \\
S40 & [1,1,4,19,100,562,3304,20071,...] & R \circ \mbox{INVERT}^3 & \alpha \\
S41 & [1,1,1,2,4,9,21,51,127,323,835,...] & R^2 \circ \mbox{INVERT} & \alpha \\
S42 & [1,1,1,1,2,4,8,17,37,82,185,423,...] & R^3 \circ \mbox{INVERT} & \alpha \\
S43 & [1,1,1,1,1,2,4,8,16,33,69,146,...] & R^4 \circ \mbox{INVERT} & \alpha \\
S39 & [1,1,3,11,45,197,903,4279,20793,...] & M^{-1} \circ \mbox{INVERT} & \beta \\
\\
S44 & [1,2,4,7,10,12,18,40,44,45,...] & \mbox{REVERT} & \epsilon \\
\\
S15 & [1,1,2,6,24,120,720,5040,...] & R \circ \mbox{EXP}  & \alpha \\
S45 & [1,1,1,2,5,16,61,272,1385,...] & R^2 \circ \mbox{EXP} & \alpha \\
S46 & [1,1,1,1,2,5,15,53,213,961,...] & R^3 \circ \mbox{EXP} & \alpha \\
S47 & [1,1,3,14,89,716,6967,79524,...] & R \circ \mbox{EXP}^2 & \alpha \\
S48 & [1,1,4,26,236,2752,39208,...] & M^{-1} \circ \mbox{EXP} & \beta \\
\end{array}
$$
\end{table}

\begin{table}[htb]
\begin{center}
{\bf Table~II~~~Fractional eigensequences} \\
\end{center}

$$
\begin{array}{lllc}
\multicolumn{1}{c}{\#} & \multicolumn{1}{c}{\mbox{Sequence}} & \multicolumn{1}{c}{\mbox{Operator}} & \multicolumn{1}{c}{\mbox{Property}} \\ [+.2in] \\
S49 & \left[ 1, - \frac{1}{2}, 0, \frac{1}{4} , 0, - \frac{1}{2}, 0, \frac{17}{8}, 0, - \frac{31}{2}, 0, \frac{691}{4}, ... \right] & N^{-1} \circ \mbox{BINOMIAL} & \beta \\ [+.1in]
S50 & \left[ 1, - \frac{1}{2}, \frac{1}{4} , \frac{1}{2} , - \frac{19}{8} , - \frac{39}{16}, \frac{2623}{32}, - \frac{365}{4}, ... \right] & N^{-1} \circ \mbox{STIRLING} & \beta \\ [+.1in]
S51 & \left[ 1, \frac{1}{2}, \frac{1}{2}, \frac{1}{4}, \frac{1}{2}, 0, \frac{1}{2}, \frac{1}{8} , \frac{1}{4}, 0, \frac{1}{2}, - \frac{1}{8}, \frac{1}{2}, 0, ... \right] & N^{-1} \circ \mbox{M\"{O}BIUS} & \beta \\ [+.1in]

S52 & \left[ 1, - \frac{1}{2}, \frac{1}{4} , \frac{1}{8} , - \frac{13}{16}, \frac{47}{32}, \frac{73}{64}, - \frac{2447}{128}, \frac{16811}{256} , ... \right] & N^{-1} \circ \mbox{EXP} & \beta
\end{array}
$$
\end{table}

We now comment on some of the individual sequences, especially those that have already
appeared in the literature.
\paragraph{BINOMIAL.}
S1 is the familiar sequence of Bell numbers,
with $\sA (x) = e^{e^x -1}$, $\sL (x) = e^x$
\cite{Bell34}, \cite{Touc56}, \cite{Rota64}, \cite[210]{C1}.
$a_n$ is the number of ways of partitioning a set of $n$ objects
(Motzkin's ``sets of sets'' \cite{Motz71}).
This sequence has property $\alpha$, as indicated in the table.
In particular, $\mbox{BINOMIAL} \circ a = L \circ a$.
In fact $b= L \circ a = [1,2,5,15,52, \dd ]$
is the unique sequence $b$ that begins $1, \dd$ and has the property that the leading
diagonal of its difference table is $R \circ b$.
Although this property of the Bell numbers must be well-known, we have not seen it
mentioned before.
$a_n$ is also the number of orbits on $n$-tuples of an $n$-fold transitive group.
Sequences S2, S3, S4 (and further sequences of the same type)
appear in Kerber \cite{Kerb78} (see also Tsuzuku \cite{Tsuz61}, Foulkes \cite{Foul70}) in connection with the characterization of multiply transitive permutation groups by cycle
structure.
These sequences have $\sL (x) = e^{2x}$, $e^{3x}$, $e^{4x}$, etc., a pleasant
property not recorded in \cite{Kerb78}.
Furthermore $b= L \circ S2$ is the unique sequence $b$ that begins $1, \dd$ and
has the property that the leading diagonal of its difference table of depth 2 is $R \circ b$.
Similar properties hold for the other sequences.
Sequence S6 has $\sA (x) =1/(2-e^x)$ (from Eq.~\eqn{eq2}), and gives the
number of unlabeled planar rooted trees with $n+1$ terminal nodes all at the same height,
in which at every other height there is a node with at least two successors
\cite{Cayl91}, \cite{MoFr84}.
E.g. $a_2 =3$ refers to \\

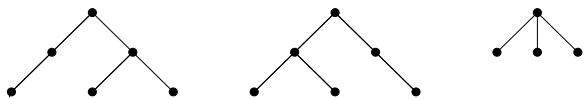
\centerline{\psfig{file=eig1.ps,width=3in}}

\vspace*{+.1in}
\noindent
Equivalently, $a_n$ is
the number of preferential
arrangements or total preorders of $n$ things (Motzkin's ``lists of sets''),
i.e. the number of linear orderings of $n$ things
using $<$ and $=$ \cite{Gros62}, \cite{Motz71}.
The orderings corresponding to the above trees are
$1<2$, $2<1$, $1=2$.
Sequence S49 has $\sA (x) =2/(1+e^x)$ and is related to the Bernoulli
numbers $B_n$ \cite[810]{AS1},
the tangent numbers $T_n$ \cite{KnBu67}, and
the $H_n$ numbers of Terrill and Terrill \cite{TeTe45} by
\begin{eqnarray*}
a_n & = & - \frac{2}{n+1} (2^{n+1} -1) B_{n+1} ~~(n \ge 0) ~, \\
& = & (-1)^{(n+1)/2} \frac{T_n}{n} ~=~
- \frac{H_{(n+1)/2}}{2^{(n+1)/2}} ~~~ (n \ge 1 ) ~.
\end{eqnarray*}
Donaghey \cite{Dona76} has studied a different kind of self-inverse sequence
associated with this transform, namely
sequences satisfying $a_n = \sum\limits_{k=0}^n (-1)^k {\binom{n}{k}} d^{n-k} a_k$,
where $d= 2 \delta$ is a constant.
\paragraph{STIRLING.}
Only one of these sequences seems to have been studied before.
For sequence S10 we have
$\sA (e^x -1) = 2 \sA (x) -1$,
$a_n = \sum\limits_{k=0}^{n-1} S(n,k) a_k$ $(n \ge 1)$.
This sequences gives the number of rooted (but not planar) trees with $n+1$
labeled nodes all at the same height
\cite{Leng84}.
E.g. $a_2 =4$ refers to \\

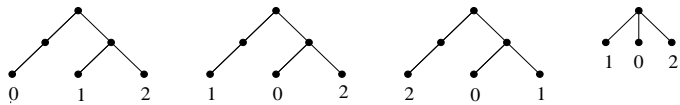
\centerline{\psfig{file=eig2.ps,width=3.5in}}

\vspace*{+.1in}
\noindent
Equivalently, $a_n$ counts a certain class of preferential
arrangements on unordered pairs of elements of an $(n+1)$-set, called
``ultradissimilarity relations'' \cite{Scha80}.

The e.g.f.'s for S7 and S50 satisfy
\begin{eqnarray*}
\sA ' (x) - \sA (e^x -1) & = & 1 ~, \\
\sA (x) + \sA (e^x -1) & = & 2
\end{eqnarray*}
respectively.
Recurrences for all these sequences follow from \eqn{eq21}--\eqn{eq23}.

The identity
$$\sum_{l=k}^{n-1}
{\binom{n-1}{l}} S_{l,k} = S_{n,k+1}$$
(\cite[p.~209]{C1})
shows that
$$\mbox{STIRLING} \circ R = R \circ \mbox{BINOMIAL} \circ \mbox{STIRLING} ~,$$
(with a little care about the zeroth term),
from which it follows that STIRLING maps the all-ones sequence to the Bell numbers (S1).
\paragraph{CONV.}
S11 is the ubiquitous sequence of Catalan numbers,
one realization of which is the number of unlabeled planar rooted trees in which
each node has $\le 2$ successors
\cite[p.~53]{C1},
\cite{Robe84}, \cite{GKP}.
Sequences S11, S12, S13, etc. satisfy
$$\mbox{CONV}^r \circ a = L \circ a ~,$$
hence $xA(x)^r = A(x) -1$, and are given by
$$a_n = \frac{1}{(r-1)n+1} {\binom{rn}{n}} ~.$$
These sequences have a long history \cite[Section~7.5, Example~5]{GKP},
\cite{Knut92}.
They also appear in the context of enumerating clusters of polygons in \cite{HPR75}.

The sequences satisfying
$$\mbox{CONV} \circ a = L^s \circ a ,~~~(s=2,3, \dd )$$
appear to be new.
The first of them, S14, has
$$A(x) = \frac{1- \sqrt{1-4x^2 -4x^3}}{2x^2} ~.$$
Except for S15, the factorial numbers, all the sequences emerging from the other
convolutional transforms are new.
\paragraph{M\"{O}BIUS.}
From here on the subscripts of the sequences start at 1.

S23, with $a_1 =1$,
$$a_{n+1} = \sum_{d|n} a_d , ~~~n \ge 1 ~,$$
is the number of planted achiral trees
\cite{GHR82}, \cite{HaRo75}.

S26, with $a_1 =1$ and
$$a_n = \sum_{d|n \atop d < n} a_d , ~~~(n \ge 2) ~,$$
is the number of ordered factorizations of $n$.
Its Dirichlet generating function
$\sD(s) = \sum\limits_{n=1}^\infty a_n / n^s$ is $2/(1+ \zeta (s)^{-1})$.
$a_n$ is also the number of ``perfect partitions'' of $n-1$, that is,
partitions of $n-1$ which contain precisely one partition of every smaller
number (i.e. a set of weights
such that every weight of $m$ grams,
$m \le n-1$, can be realized in exactly one way)
\cite{MacM91}, \cite[p.~123]{R1}, \cite[p.~126]{C1}.

The other sequences appear to be new.
Sequence S51 has Dirichlet generating function
$1/(2- \zeta (s)^{-1})$.
\paragraph{WEIGH.}
S27 gives the number of unlabeled rooted trees with no symmetries (\cite{HRS75}, see in particular Eq.~(1$^\ast$)).
That this sequence shifts left under WEIGH was already pointed out by
Cameron \cite{Came89}.
\paragraph{EULER.}
S30 counts unlabeled rooted trees, and satisfies
$$\sum_{n=0}^\infty a_{n+1} x^n = \prod_{n=1}^\infty \frac{1}{(1-x^n)^{a_n}}$$
\cite[p.~127]{R1}, \cite{HaRo75}.
S33 counts series-reduced planted trees, and $M \circ S33$ counts
series-parallel networks \cite{RiSh42},
\cite[p.~142]{R1}, \cite{Lomn72}.
The connections of these three sequences with the EULER transform have been investigated by Cameron \cite{Came87}, \cite{Came89}.
\paragraph{PARTITION.}
Every sequence converges to one of these five.
S34 is fixed,
while S35 swaps with S36, and S37 with S38.
There are no other sequences of finite order under this transformation.
In particular,
S34 is the unique sequence $[a_1, a_2, \dd ]$ with the property
that the number of partitions of $n$ into parts of sizes $\{a_1 , a_2 , \dd \}$ is equal to $a_n$.
We find it remarkable that, in spite of all the papers that have been written about partitions, these five sequences seem never to have been published before.
\paragraph{INVERT.}
Cameron \cite{Came89} has already observed that the Catalan numbers (S11) shift left under
INVERT, and that S39, the solution to Schr\"{o}der's second problem satisfies
\linebreak
$\mbox{INVERT} \circ a = M \circ a$.
Like the Catalan numbers, the latter sequence occurs in many guises, sometimes
being called super-Catalan numbers
\cite{Schr70}, \cite{Ethe40},
\cite[p.~168]{RCI},
\cite{Krew73},
\cite[p.~57]{C1},
\cite{Roge77},
\cite{Lew79},
\cite{StWa79},
\cite{GSWW}.
$a_n$ is among other things the number of rooted planar trees with $n+1$ nodes.
The o.g.f. satisfies $2A(x)^2 - (1+x)A(x) + x=0$, which implies the
inversion property mentioned above, as well as the fact that this sequence
also satisfies $\mbox{INVERT}^2 \circ a = L \circ a$.

S41 is the sequence of Motzkin numbers, again a widely-studied sequence
\cite{Motz48}, \cite{Carl69}, \cite{DoSh77},
\cite{Dona77}, \cite{Dona80}.
The o.g.f. satisfies $A(x)^2 = (1+x) (A(x) -x)$.

Sequences S42, S43 and similar sequences are related to the Motzkin numbers,
and have been studied in \cite{Wate78}, \cite{StWa79}, \cite{HSW80},
in part because they enumerate secondary structures of nucleic acids.
\paragraph{REVERT.}
S44 is new, and no combinatorial setting is presently known for it.
\paragraph{EXP.}
S45 is the sequence of Euler numbers,
with e.g.f. $(1+ \sin x)/\cos x$ \cite[p.~110]{NET},
\cite[p.~262]{DKB},
\cite[p.~259]{C1}.

S48 is the solution to Schr\"{o}der's fourth problem, which is also
connected with series-parallel networks
\cite[p.~197]{RCI}, \cite[p.~224]{C1},
\cite{Lomn72}.
The e.g.f. satisfies $\exp \sA (x) = 2 \sA (x) +1 -x$
\cite[p.~224, Eq.~(aa)]{C1}.

S52 arises in an unexpected context.
If we multiply $a_n$ by $2^{n-1}$ we obtain an integer sequence
$[1,-1,1,1,-13,47,73, \dd ]$, sequence 2093 in \cite{HIS}.
This occurs in numerical analysis in the computation of Airey's converging
factor for asymptotic series, both for the Weber parabolic cylinder function and for
the exponential integral \cite{Aire37}, \cite{Mill52},
\cite{Wynn63}, \cite{Murn72}.
The simple equation for the e.g.f. for S52 which follows from the
eigen-sequence property,
$$\exp \sA (x) = 1+ 2x - \sA (x) ~,$$
does not seem to be given in any of these references.
Once again the other sequences are new.
\subsection*{Acknowledgements}
\hsp
We thank Peter Cameron for some valuable comments on a preliminary version of
the manuscript.

\end{document}